\documentclass[a4paper,14pt]{article}

\usepackage{multicol}
\usepackage{amsthm}

\usepackage[dvips]{graphicx,color}
\usepackage{amssymb}
\usepackage{amsmath}

\begin{document}

\fontsize{14pt}{16.5pt}\selectfont

\begin{center}
\bf{A product space of $\{0$, $1\}$ and an abstract polycrystal}
\end{center}

\fontsize{12pt}{15pt}\selectfont
\begin{center}
Akihiko Kitada$^{1}$, Shousuke Ohmori$^{2}$, Tomoyuki Yamamoto$^{1,2,*}$\\ 
\end{center}

\noindent
$^1$\it{Institute of Condensed-Matter Science, Comprehensive Resaerch Organization, Waseda University,
3-4-1 Okubo, Shinjuku-ku, Tokyo 169-8555, Japan}\\
$^2$\it{Faculty of Science and Engineering, Waseda~University, 3-4-1 Okubo, Shinjuku-ku, Tokyo 169-8555, Japan}\\
*corresponding author: tymmt@waseda.jp\\
~~\\

\rm
\noindent
{\bf Abstract}\\
A partition of a $\Lambda (Card\Lambda \succ \aleph _0)$-product space of $\{0,1\}$ defines an abstract polycrystal composed of abstract singlecrystals a decomposition space (space of equivalence classes) of each of which is self-similar.
%\end{abstract}

\section{Introduction}

Let $(X, \tau)=(\{0,1\}^\Lambda ,\tau_0^\Lambda ),~Card \Lambda \succ \aleph _0$
(aleph zero) be the $\Lambda -$product space of $(\{0,1\},\tau_0)$ where $\tau_0$ is a 
discrete topology for $\{0,1\}$. The topological space $(X,\tau)$ needs not to be metrizable, 
that is, the relation $Card \Lambda \succ \aleph $ (aleph) may hold. In the present report, assuming the lattice point of crystal to be a map $x:\Lambda \rightarrow \{0,1\}$ (for example, $\Lambda =\bf{N}, \bf{R},\dots$), we 
will mathematically confirm the existence of a partition $\{X_1,\dots,X_n\}$ of $X$, $X_i \in 
{(\tau\cap \Im)}$\footnote{The notaion $\tau\cap \Im$ denotes the set of all closed and open sets (clopen sets) of $(X,\tau)$.}$-\{\phi\}$, a decomposition space (i.e, a space of equivalence 
classes) of each element $X_i$ of which is self-similar. The partition of $X$ can be regarded as a kind of polycrystal\footnote{We may replace a polycrystal with a tiling and a singlecrystal with a tile, respectively.} in an abstract sense each abstract singlecrystal $X_i=(X_i,\tau_{X_i})$\footnote{$\tau_A=\{u\cap A;u\in \tau\}$ for $A\subset X$.} of which is characterized by its self-similar decomposition space. Namely, all points in $X_i$ are classified into equivalence classes and then the equivalence classes coalesce to form a self-similar structure.

Some easily verified statements are summarized in the next section for the preliminaries.

\section{Preliminaries}
\begin{description}
\item[2-a)] Any zero-dimensional (0-dim), perfect, T$_0$(necessarily T$_2$)-space $(X,\tau)$ has its partition $\{X_1,\dots,X_n\}, X_i\in(\tau\cap \Im)-\{\phi\}$ for any $n$. Here $X_i\cap X_{i'}=\phi$ for $i\not =i'$ and $\displaystyle\bigcup _{i\in \bar n}X_i=X$\footnote{$\bar n=\{1,\dots,n\}$.}. Each subspace $(X_i,\tau_{X_i})$ is a 0-dim, perfect, T$_2$-space.
\item[2-b)] Let $X$ be a 0-dim, perfect, compact T$_2$-space. Then, for any compact metric space $Y$, there exists a continuous map $f$ from $X$ onto $Y$.
\item[2-c)] If $f:(X,\tau)\rightarrow (Y,\tau')$ is a quotient map, then the map $h:(Y,\tau')\rightarrow (\mathcal{D}_f,\tau(\mathcal{D}_f)), y\mapsto f^{-1}(y)$ is a homeomorphism. Here, the decomposition $\mathcal{D}_f$ of $X$ and the decomposition topology $\tau(\mathcal{D}_f)$ are given by $\mathcal{D}_f=\{f^{-1}(y)\subset X;y\in Y\}$ and $\tau(\mathcal{D}_f)=\{\mathcal{U}\subset \mathcal{D}_f;\bigcup \mathcal{U}\in\tau \}$, respectively.
\item[2-d)] A space which is homeomorphic to a self-similar space is self-similar.
\item[2-e)] Let $(X,\tau)$ be a 0-dim, perfect, compact, not-metrizable T$_2$-space. From the Urysohn's metrization theorem $(X,\tau)$ is not second countable. In the partition $\{(X_1,\tau_{X_1}),\dots, \\ (X_n,\tau_{X_n})\}$ of $(X,\tau)$, there exists a number $i_0\in \bar n$ such that $(X_{i_0},\tau_{X_{i_0}})$ is not second countable. The subspace $(X_{i_0},\tau_{X_{i_0}})$ is a 0-dim, perfect, compact not-metrizable T$_2$-space.
\end{description}

\section{A partition of the space $(\{0,1\}^\Lambda ,\tau_0^\Lambda )$}

Since $(X,\tau)=(\{0,1\}^\Lambda ,\tau_0^\Lambda ),\tau_0=2^{\{0,1\}},Card 
\Lambda \succ \aleph _0$ is easily verified to be a 0-dim, perfect, compact T$_2$-space, from 
{\bf 2-a)} there exists a partition $\{(X_1,\tau_{X_1}),\dots,(X_n,\tau_{X_n})\}$ of $(X,\tau)$ where each $(X_i,\tau_{X_i})$ is a 0-dim, perfect, 
compact T$_2$-space.\footnote{If $Card \Lambda \succ \aleph $, as mentioned in {\bf 
2-e)}, there exists 0-dim, perfect,compact, {\it not-metrizable} T$_2$-space 
$(X_{i_0},\tau_{X_{i_0}})$.} Then, from {\bf 2-b)}, there exists a continuous map $f_i$ 
from $(X_i,\tau_{X_i})$ onto any compact, self-similar metric space $(Y,\tau_d)$. Since $(X_i,\tau_{X_i})$ is a compact space and $Y$ is a T$_2$-space, the map $f:(X_i,\tau_{X_i})\rightarrow (Y,\tau_d)$ is a quotient map. Therefore, from {\bf 2-c)}, the map $h:(Y,\tau_d)\rightarrow (\mathcal{D}_{f_i},\tau(\mathcal{D}_{f_i})),y\mapsto f_i^{-1}(y)$ must be a homeomorphism. Since $(Y,\tau_d)$ is self-similar, according to {\bf 2-d)}, the decomposition space $(\mathcal{D}_{f_i},\tau(\mathcal{D}_{f_i}))$ of $(X_i,\tau_{X_i})$ is also self-similar. Here, we note that the decomposition space $\mathcal{D}_{f_i}$ of $X_i$ is not a trivial one $\{\{x\};x\in X\}$ especially for a self-similar, connected space such as the Sierpi\'nski carpet (S.B.Nadler Jr, 1992). In fact, disconnected space $\{\{x\};x\in X\}$ (a decomposition space of $X_i$) is never homeomorphic to a connected space.

Regarding each subspace $(X_i,\tau_{X_i}),i\in \bar n$ an abstract singlecrystal we obtain an abstract polycrystal $X$ composed of $(X_i,\tau_{X_i})$ which has a self-similar decomposition space.

Finally, we note that in the above discussions we can replace the self-similar space with a compact substance in the materials science such as dendrite \cite{dendrite} and then, we can obtain an abstract polycrystal $(\{0,1\}^\Lambda ,\tau_0^\Lambda )$ composed of abstract singlecrystal whose decomposition space is characterized by a dendrite.

\section*{Acknowledgement}
The authors are grateful to Dr. Y. Yamashita and prof. H. Nagahama at Tohoku university for useful suggestions and encouragements.\\
 
%\section{references}

\end{document}